\documentclass[12pt,a4paper]{amsart}
\usepackage{amssymb,amsmath}

\headsep=1cm \numberwithin{equation}{section}
\newtheorem{theorem}{Theorem}[section]
\newtheorem{lemma}[theorem]{Lemma}

\newtheorem{corollary}[theorem]{Corollary}

\theoremstyle{definition}

\theoremstyle{remark}

\newcommand{\Ass}{\operatorname{Ass}}

\newcommand{\Ht}{\operatorname{ht}}

\newcommand{\Rad}{\operatorname{Rad}}

\newcommand{\fm}{\frak{m}}
\newcommand{\fp}{\frak{p}}

\begin{document}

\author[Divaani-Aazar and Tousi ]{Kamran Divaani-Aazar and Massoud Tousi }
\title[Localization of tight closure in two-dimensional rings]
{Localization of tight closure in two-dimensional rings}

\address{K. Divaani-Aazar, Department of Mathematics, Az-Zahra University,
Vanak, Post Code 19834, Tehran, Iran and Institute for Studies in
Theoretical Physics and Mathematics, P. O. Box 19395-5746, Tehran,
Iran.} \email{kdivaani@ipm.ir}

\address{M. Tousi, Institute for Studies in Theoretical Physics and
Mathematics, P. O. Box 19395-5746, Tehran, Iran and Department of
Mathematics, Shahid Beheshti University, Tehran, Iran.}
\email{mtousi@ipm.ir}

 \maketitle

\begin{abstract} It is shown that tight closure commutes with
 localization in any two dimensional ring $R$ of prime characteristic
 if either $R$ is a Nagata ring or $R$ possesses a weak test
 element. Moreover, it is proved that tight closure commutes with
 localization at height one prime ideals in any ring of prime
 characteristic.
\end{abstract}

\section{Introduction}

 Throughout this paper, $R$ is a commutative Noetherian ring (with identity)
 of prime characteristic $p$. The theory of tight closure was
 introduced by Hochster and Huneke [{\bf 2}]. There are
 many applications for this notion in both commutative algebra
 and algebraic geometry. However, there are many basic open
 questions concerning tight closure. One of the essential
 questions is whether tight closure commutes with localization.
 For an expository account on tight closure, we refer the reader
 to [{\bf 3}] or [{\bf 8}].

 In the sequel, $R^{\circ}$ denotes the set of elements of $R$ which
 are not contained in any minimal prime ideal of $R$. We use the
 letter $q$ for nonnegative powers $p^e$ of $p$. Let $I$ be an
 ideal of $R$ and $I^{[q]}$ the ideal generated by $q$-th powers
 of elements of $I$. Then $I^*$, {\it tight closure} of $I$ is the set of
 all elements $x\in R$ for which there exists $c\in R^{\circ}$ such
 that $cx^q\in I^{[q]}$ for all $q\gg 0$. Also, for a nonnegative
 power $q'$ of $p$ an element $c\in R^{\circ}$ is called $q'$-weak test
 element, if for any ideal $I$ of $R$ and any element $x\in I^*$,
 we have $cx^q\in I^{[q]}$ for all $q\geq q'$.

 We say that tight closure commutes with  localization for the
 ideal $I$, if for any multiplicative system $W$ in $R$,
 $I^*R_W=(IR_W)^*$. It is conjectured that tight closure commutes
 localization in general situation. There are some related conjectures that a
 positive answer to each of them will yield solution to the localization
 problem. For example, it is an open question that for any ideal $I$ of a domain $R$,
 $I^*=IR^+\cap R$, where $R^+$ denote the integral
 closure of $R$ in an algebraic closure of its fraction field.
 An positive solution to this problem implies a solution
 to localization problem (see [{\bf 3}]). Also [{\bf 7}], if the Frobenius
 powers of the proper ideal $I$ of $R$ have {\it linear growth of primary
 decompositions}, then tight closure of $I$ commutes with localization at a
 multiplicative system consisting of the powers of a single element of
 $R$.

 Tight closure commutes with localization in several important special cases.
 For example, it is known that tight closure commutes with localizations on
 principal ideals and also on ideals generated by regular sequences (see e.g. [{\bf 3}]).
 We refer the reader to [{\bf 8}] for a survey of results on localization
 problem. Also, Smith [{\bf 6}] proved that tight closure commutes with
 localization in affine rings which are quotients of a polynomial ring over a
 field, by a binomial ideal. She proved this result, by showing that if for
 any minimal prime ideals $\fp$ of the ring $R$, the quotient $R/\fp$ has a finite extension
 domain in which tight closure commutes with localization, then tight closure
 commutes with localization in $R$ itself. This fact is essential in the proof of
 our main result:

\begin{theorem} Let $R$ be a Noetherian ring of prime
 characteristic. Then tight closure commutes with localization in
 $R$ in the following cases: \\
i) $\dim R < 2$. \\
ii) R is a two dimensional local ring. \\
iii) $\dim R=2$ and either $R$ is a Nagata ring or $R$ possesses a
$q'$-weak test element.\\
\end{theorem}
Note that, by [{\bf 4}, Theorem 78 and Definition 34.A], any
excellent ring is Nagata. \\
It seems that the solution of
localization problem in dimension two has been basically known (at
least in the case of excellent rings) by most experts that have
attacked the problem. However, since the solution did not appear
in any article, at least according to our knowledge, the main
achievement of this paper is that it fills a gap in the
literature.

\section{The proof}

To prove the theorem, we proceed through the following lemmas,
some of which may be of independent interest.

\begin{lemma} Let $I$ be  an ideal of $R$ and $W$ a multiplicative
system in $R$. Suppose there exists $w\in W$ such that $w^q(
I^{[q]}R_W\cap R)\subseteq  I^{[q]}$ for all $q\gg 0$. Then
$I^*R_W=(IR_W)^*$.
\end{lemma}

{\bf Proof.} Clearly $I^*R_W\subseteq (IR_W)^*$. Now, let $x/1\in
(IR_W)^*$. Then there is $c\in (R_W)^{\circ}$ such that
$c(x/1)^q\in (IR_W)^{[q]}$ for all $q\gg 0$. It is easy to see
that in fact we can chose $c$ in $R^{\circ}$.  Hence $cx^q\in
I^{[q]}R_W\cap R$ and so $w^qcx^q=c(wx)^q\in I^{[q]}$ for all
$q\gg 0$. Thus $wx\in I^*$, and therefore $x\in I^*R_W$, as
required. $\blacksquare$

By [{\bf 2}, Proposition 4.14], for any maximal ideal and any
$\fm$-primary ideal $I$ of $R$, $I^*R_{\fm}=(IR_{\fm})^*$. The
following extend this fact.

\begin{lemma} Let $I$ be an ideal of $R$. Let $\fp$ be a maximal
ideal of $R$ which is minimal over $I$. Then $I
^*R_{\fp}=(IR_{\fp})^*$. In particular, if $R/I$ is an Artinian
ring, then tight closure commutes with localization for $I$.
\end{lemma}

{\bf Proof.} Since $IR_{\fp}$ is $\fp R_{\fp}$-primary,  it
follows that $IR_{\fp}$ contains some power of  $\fp R_{\fp}$.
Suppose $\fp ^kR_{\fp}\subseteq IR_{\fp}$. Then $s\fp ^k\subseteq
I$ for some $s\in R\smallsetminus\fp$. Assume that $I$ is
generated by $t$ elements and let $w=s^t$ and $n=tk$. Then $w^q$
multiplies the $nq$-th power of $\fp$ into $I^{[q]}$ for all $q$.
By Lemma 2.1, it suffices to show that $$w^q( I^{[q]}R_{\fp}\cap
R)\subseteq I^{[q]},$$ for all $q$. To this end, let $x\in
(I^{[q]}R_{\fp}\cap R)$ for some $q$. Then there exists $s_q\in
R\smallsetminus\fp$ such that $s_qx\in  I^{[q]}$. Since $Rs_q+
\fp^{nq}=R$, there are $\alpha_q$ in $R$ and $\beta_q$ in
$\fp^{nq}$ such that $\alpha_qs_q+ \beta_q=1$. Then
$$w^qx=\alpha_qs_qw^qx+ \beta_qw^qx\in I^{[q]}.$$ Therefore $w^q(
I^{[q]}R_{\fp}\cap R)\subseteq  I^{[q]}$.

For the second assertion, first note that, by [{\bf 1}, Lemma
3.5(a)], it is enough to treat only the localization at prime
ideals of $R$. Also, note that for  any prime ideal $I$ which does
not contain $I$, we have $I^*R_{\fp}=(IR_{\fp})^*=R_{\fp}$. Hence
the claim follows by the first part of the lemma. $\blacksquare$

It follows from Lemma 2.2 that tight closure commutes with
localization in any domain of dimension less  that or equal to
one. Thus by using [{\bf 6}, Lemma 1], we deduce the following
result.

\begin{corollary} Assume $R$ is a ring with $\dim R\leq 1$. Then
tight closure commutes with localization in $R$.
\end{corollary}

\begin{lemma} Let $I$ be an ideal of the integral domain $R$.
Then $I ^*R_{\fp}=(IR_{\fp})^*$ for any height one prime ideal
$\fp$ of $R$.
\end{lemma}

{\bf Proof.} Let $\fp$ be a height one prime ideal of $R$. Let
$R'$ be the integral closure of $R$ in its field of fractions. It
is known that the normalization of any one dimensional Noetherian
domain is Noetherian. Since $(R')_{\fp}$ is the integral closure
of the domain $R_{\fp}$ in its field of fractions, it follows that
$(R')_{\fp}$ is a Noetherian normal domain of dimension one. Hence
$(R')_{\fp}$ is  a regular ring. This implies that every ideal of
$(R')_{\fp}$ is tightly closed.

Now, let $x/1\in (IR_{\fp})^*$. Then there is a non-zero element
$c$ in $R$ such that $c(x/1)^q\in (IR_{\fp})^{[q]}$ for all $q$.
But for each $q$, $(IR_{\fp})^{[q]}\subseteq (I(R')_{\fp})^{[q]}$,
and so $x/1\in (I(R')_{\fp})^*=I(R')_{\fp}$. Hence there is $s\in
R\smallsetminus\fp$ such that $sx\in IR'$. This implies that
$x/1\in I^*R_{\fp}$, because $IR'\cap R\subseteq I^*$, by [{\bf
3}, Page 15]. $\blacksquare$

Lemma 2.4 has the following interesting conclusion.

\begin{corollary} Let $R$ be a ring and $I$ an ideal of $R$. For any
height one prime ideal $\fp$ of $R$, we have
$I^*R_{\fp}=(IR_{\fp})^*$.
\end{corollary}

{\bf Proof.} Let $\{\fp_1,\fp_2,\dots ,\fp_n \}$ be the set of
minimal prime ideals of $R$, which are contained in $\fp$. Fix
$1\leq i\leq n$. By Lemma 2.4, we have
$$(I+\fp_i/\fp_i)^*(R/\fp_i)_{\fp/\fp_i}=((I+\fp_i/\fp_i)(R/\fp_i)
_{\fp/\fp_i})^*.$$
Now, by following the argument of [{\bf 6}, Lemma 1], it turns out
that $I^*R_{\fp}=(IR_{\fp})^*$. $\blacksquare$

The following is the only technical tool remaining in order to
prove Theorem 1.1. Some tricks in the argument of the following
result is very close to those which are used in [{\bf 9},
Proposition 1.2].

\begin{lemma} Let $R$ be a two dimensional normal ring. Let $I$ be
a height one ideal of $R$ and $\fp$ a height two prime ideal of
$R$ containing $I$. Then $I^*R_{\fp}=(IR_{\fp})^*$.
\end{lemma}

{\bf Proof.} If $\fp$ is minimal over $I$, then the assertion
follows by Lemma 2.2. Thus in the sequel, we assume that $\fp$ is
not minimal over $I$. Suppose that $\{\fp_1,\fp_2,\dots ,\fp_n \}$
is the set of minimal prime ideals of $I$, which are contained in
$\fp$. Let $I^{[q]}=\cap_{i=1}^{t_q}Q_{iq}$ be a minimal primary
decomposition of the ideal $I^{[q]}$, with
$\Rad(Q_{iq})=\fp_{iq}$. Then $$I^{[q]}R_{\fp}\cap
R=\bigcap_{\fp_{iq}\subseteq \fp}Q_{iq}\subseteq
\bigcap_{i=1}^n(I^{[q]}R_{\fp_i}\cap R).$$

Since $R$ is reduced, it follows that each associated prime ideal
of $R$ is minimal. Hence, by Prime Avoidance Theorem, we can
deduce that $I$ can be generated by regular elements in $R$.
Because each $R_{\fp_i}$ is a DVR, there are $a_1,a_2,\dots ,a_n$
in $I$ such that each $a_i$ is a regular in $R$ and
$IR_{\fp_i}=a_iR_{\fp_i}$. It follows that $\fp_i$ is a minimal
over $a_iR$. Let $\fp_{i2},\fp_{i3}\dots ,\fp_{in_i}$ be the other
associated prime ideals of the ideal $a_iR$. Take $z_i\in
\bigcap_{j=2}^{n_i}\fp_{ij}\smallsetminus\fp_i$. For an ideal $J$
of $R$, we denote $\bigcup_{k\in \mathbb{N}} J:z_i^k$, by
$J:<z_i>$. Since $a_i$ is a regular element of $R$, it follows
that for each $q$, $\Ass_R(R/a_i^qR)=\Ass_R(R/a_iR)$, and so one
can check easily that $$a_i^qR_{\fp_i}\cap R=a_i^qR:<z_i>.$$

By [{\bf 3}, Exercise 4.2], there exists an integer $c_i$ such
that $a_i^qR:<z_i>=a_i^qR:z_i^{c_iq}$ for all $q$. Therefore
$$I^{[q]}R_{\fp}\cap R\subseteq \bigcap_{i=1}^n(a_i^qR:z_i^{c_iq}),$$
for all $q$. Because the ideal $\sum_{i=1}^nz_i^{c_i}R$ is not
contained in the union of $\fp_i$'s, $i=1,2,\dots ,n$, there are
elements $r_1,r_2,\dots ,r_n$ in $R$ such that
$$\alpha=r_1z_1^{c_1}+r_2z_2^{c_2}+\dots +r_nz_n^{c_n} \notin
\bigcup_{i=1}^n \fp_i.$$ Let $x\in I^{[q]}R_{\fp}\cap R$ for some
$q$. Then for each $i=1,2,\dots ,n$, we have $z_i^{c_iq}x\in
a_i^qR\subseteq I^{[q]}$. Thus $\alpha^qx\in I^{[q]}$. If $\alpha
\notin \fp$, then by Lemma 2.1, it turns out that
$I^*R_{\fp}=(IR_{\fp})^*$.

Now, assume that $\alpha \in \fp$. Since $\alpha \notin
\cup_{i=1}^n \fp_i$, it turns out that  $\fp$ is minimal over
$\alpha R+I$. Using the similar argument as in the proof of Lemma
2.2, we can deduce that there are $s\in R\smallsetminus\fp$ and
$l\in \mathbb{N}$ such that $s^q\fp^{lq}\subseteq
\alpha^qR+I^{[q]}$ for all $q$. Since $x\in I^{[q]}R_{\fp}\cap R$,
there is $w_q\in R\smallsetminus\fp$ such that $w_qx\in I^{[q]}$.
Because $\fp^{lq}+Rw_q=R$, there are $\beta_q\in \fp^{lq}$ and
$\gamma_q\in R$ such that $1=\beta_q+\gamma_qw_q$, and so
$$s^qx=\beta_qs^qx+s^q \gamma_qw_qx\in I^{[q]}.$$
Hence $s^q(I^{[q]}R_{\fp}\cap R)\subseteq I^{[q]}$ for all $q$ and
so the claim follows by using Lemma 2.1 again. $\blacksquare$

Recall that an element $c\in R^{\circ}$ is called $q'$-weak test
element, if for any ideal $I$ of $R$ and any element $x\in I^*$,
we have $cx^q\in I^{[q]}$ for all $q\geq q'$. The following
improves [{\bf 8}, Proposition 6.5].

\begin{lemma} Let $R\subseteq S$ be an integral extension of
Noetherian domains. Suppose that $R$ possesses a $q'$-weak test
element. Then for any ideal $I$ of R, $(IS)^*\cap R\subseteq I^*$.
\end{lemma}

{\bf Proof.} Let $x\in (IS)^*\cap R$. Then there is a non-zero
element $d$ of $S$ such that $dx^q\in (IS)^{[q]}$ for all $q\gg
0$. Since $S$ is integral over $R$, there are  $a_0,a_1,\dots ,a_n
\in R$ such that $$d^n+a_1d^{n-1}+\dots +a_0=0.$$ We may and do
assume that $a_0$ is non-zero. Then $a_0x^q\in I^{[q]}S\cap R$ for
all $q\gg 0$.  Therefore, by [{\bf 3}, Page 15], $a_0x^q\in
(I^{[q]})^*$ for all $q\gg 0$. Now it follows from [{\bf 2}, Lemma
8.16], that $x\in I^*$. $\blacksquare$

Now, we are ready to conclude the Theorem 1.1.

{\it Proof of Theorem 1.1.} If $R$ has a $q'$-weak test element,
then it follows by [{\bf 1}, Lemma 2.10(a)], that for any minimal
prime ideal $\fp$, the domain $R/\fp$ has also a $q'$-weak test
element. Thus, in view of [{\bf 6}, Lemma 1], we can assume that
$R$ is a domain. Also, by [{\bf 1}, Lemma 3.5(a)], it suffices to
consider only localization at prime ideals.\\
The case (i) holds by Corollary 2.3 and the case (ii) follows by
lemma 2.4.

Now suppose that either $R$ is a Nagata ring or possesses a
$q'$-weak test element. By [{\bf 5}, Theorem 33.12], the integral
closure of a Noetherian two dimensional domain in its field of
fractions is Noetherian. Let $R'$ denote the integral closure of
$R$ in its field of fractions. Let $I$ be a non-zero ideal of $R$.
If $R$ possesses a $q'$-weak test element, then $(IR')^*\cap
R\subseteq I^*$, by Lemma 2.7. If $R$ is a Nagata ring, then $R'$
is a finite extension of $R$, and so $(IR')^*\cap R\subseteq I^*$,
by [{\bf 3}, Theorem 1.7]. Therefore in both case, by adopting the
argument of [{\bf 6}, Lemma 2], we can deduce that tight closure
commutes with localization, if the same happens in $R'$. Hence, in
the sequel, we suppose that $R$ is a Noetherian normal domain of
dimension two. If $\Ht (I)=2$, then by Lemma 2.2, tight closure
commutes with localization for $I$. Hence we may assume that $\Ht
(I)=1$. Let $\fp$ be a prime ideal of $R$. If $I$ is not contained
in $\fp$, then $I^*R_{\fp}=(IR_{\fp})^*=R_{\fp}$. Therefore the
proof is complete by Lemmas 2.4 and 2.6. $\blacksquare$

{\bf Acknowledgments.} The second named author was in part
supported by a grant from IPM. The authors thank Irena Swanson for
careful reading of an earlier version of this paper and helpful
lectures on tight closure, she presented at IPM in winter 2002. We
would also like to thank Karen Smith for pointing out that in
Corollary 2.5, there is no need to assume that $R$ is excellent.

%%%%%%%%%%%%%%%%%%%%%%%%%%%%%%%%%%%%%%%%%%%%%%%%%%%%%%%%%%%%%%%%%%%%%%%%%%%%%

\end{document}